\documentclass{article}

\usepackage[margin=0.95in]{geometry} 
\usepackage{amsmath,amsthm,amssymb,amsfonts}
\usepackage{xcolor} 
\newtheoremstyle{mytheoremstyle} 
    {\topsep}                    
    {\topsep}                    
    {}                   
    {}                           
    {\scshape}                   
    {.}                          
    {.5em}                       
    {}

\theoremstyle{mytheoremstyle}
\def\blockg#1{ %
 \bigskip 
 \fboxsep0.5em
 \noindent\fcolorbox{gray!40}{gray!20}{
 \begin{minipage}{\dimexpr\linewidth-2em}
    \noindent\quad#1
 \end{minipage}}
 \par\bigskip}

\newtheorem*{theorem}{Theorem}
\newtheorem*{definition}{Definition}
\newtheorem*{note}{Note}
\newtheorem*{corollary}{Corollary}
\newtheorem*{lemma}{Lemma}
\newtheorem*{proposition}{Proposition}

\newcommand{\R}{\mathbb{R}}

\newcommand{\Z}{\mathbb{Z}}

\newcommand{\be}{\begin{enumerate}}
\newcommand{\ee}{\end{enumerate}}
\newcommand{\bi}{\begin{itemize}}
\newcommand{\ei}{\end{itemize}}
\newcommand{\bq}{\begin{equation*}}
\newcommand{\eq}{\end{equation*}}
\newcommand{\bbm}{\begin{bmatrix}}
\newcommand{\ebm}{\end{bmatrix}}
\newcommand{\bsm}{\left[ \begin{smallmatrix}}
\newcommand{\esm}{\end{smallmatrix} \right]}

\newcommand{\ol}{\overline}

\newcommand{\bd}{\begin{definition}}
\newcommand{\ed}{\end{definition}}
\newcommand{\bn}{\begin{note}}
\newcommand{\en}{\end{note}}
\newcommand{\bl}{\begin{lemma}}
\newcommand{\el}{\end{lemma}}
\newcommand{\bt}{\begin{theorem}}
\newcommand{\et}{\end{theorem}}
\newcommand{\bp}{\noindent{\bf Proof.} }
\newcommand{\ep}{$\blacksquare$\\}
\newcommand{\bx}{\noindent{\it Example.} }
\newcommand{\ex}{$\square$\\}
\newcommand{\bc}{\begin{corollary}}
\newcommand{\ec}{\end{corollary}}
\newcommand{\bprop}{\begin{proposition}}
\newcommand{\eprop}{\end{proposition}}

\title{\textbf{Constructing an orthonormal wavelet from an MRA}}
\author{Kwok Hao Lee, Guido L. Weiss}
\begin{document}
\maketitle

\begin{abstract}
Multiresolution Analysis (MRA) wavelets have important applications in image processing and signal decomposition. In this article, we follow closely the approach in Hernandez and Weiss's seminal text [1] in describing the construction of an orthonormal wavelet from an multi-resolution analysis (MRA), a family of subspaces of $L^2(\R)$ satisfying certain properties. We assume the reader has a modest background in undergraduate analysis and measure theory.
\end{abstract}

\section{Preliminaries}

Before our work begins in earnest, it is necessary for us to give several basic definitions and prove some easy results about orthonormal systems.

\bn 
In the sections to follow, we will be using Fourier transforms of the form
\bq
\hat f(\xi) = \int_{\R} f(x) e^{-2\pi i x \xi} \, dx.
\eq
\en

\blockg{\bd
A \textbf{wavelet} is a function $\psi \in L^2(\R)$ such that
\bq
\psi_{j,k}(x) = 2^{j/2}\psi(2^j x - k), \quad j,k \in \Z,
\eq
is an orthonormal basis for $L^2(\R)$.
\ed}

\bx
The \textbf{Haar Wavelet},
\bq
\psi(x)= \left\{
\begin{array}{ll}
1 & : 0 \leq x < 1/2 \\
-1 & : 1/2 \leq x < 1 \\
0 & : \text{otherwise, }
\end{array} \right.
\eq
is an orthonormal wavelet for $L^2(\R)$; in fact, $\{\psi_{j,k}: j,k \in \Z \}$ is an orthonormal basis for $L^2(\R)$.
\ex

\blockg{\bd
A \textbf{multiresolution analysis} (MRA) comprises a sequence of closed subspaces $V_j$, $j \in \Z$, of $L^2(\R)$ satisfying
\begin{align*}
& V_j \subseteq V_{j+1} \quad \forall \, j \in \Z; \\
& f(x) \in V_j \iff f(2x) \in V_{j+1} \quad \forall \, j \in \Z; \\
& \bigcap_{j \in \Z} V_j = \{0\}; \\
& \ol{\bigcup_{j \in \Z} V_j} = L^2(\R); \\
& \exists \, \varphi \in V_0 \text{ s.t. } \{\varphi(x - k): k \in \Z \} \text{ is an orthonormal basis for } V_0. \tag{*}
\end{align*}
The function $\varphi$ whose existence is asserted in (*) is a \textbf{scaling function} of the given MRA.
\ed}

\blockg{\bl[1] Suppose $g \in L^2(\R)$. Then $\{g(\cdot - k) \, | \, k \in \Z \}$ is an orthonormal system iff
\bq
\sum_{k \in \Z} |\hat{g}(\xi - k)|^2 = 1 \quad \text{a.e.}
\eq
\el}

\bp 
Suppose $\{g(\cdot - k) \, | \, k \in \Z \}$ were an orthonormal system. Then
\begin{align*}
\delta_{k,0} &=\int_{\R} g(x) \overline{g(x-k)} \, dx = \int_{\R} |\hat{g}(\xi)|^2 e^{2\pi ik \xi} \, d \xi \quad \text{(taking inverse Fourier Transforms)}\\
&= \sum_{l \in \Z} \int_0^1 |\hat{g}(\mu - l)|^2 e^{2 \pi ik \mu} \, d \mu \\
&= \int_0^1 \sum_{l \in \Z} |\hat{g}(\mu - l)|^2 e^{2 \pi ik \mu} \, d \mu. \quad \text{(Lebesgue Dominated Convergence Theorem)}
\end{align*}
Since $\sum_{k \in \Z} |\hat g(\mu - k)|^2$ is 1-periodic, it equals $1$ a.e. since it has Fourier coefficient $1$ at frequency $k=0$ and all the other coefficients are zero. The converse follows immediately from reversing the direction of the previous argument.
\ep

Associated with Lemma 1, we have

\blockg{\bc[2] If $g \in L^2(\R)$ and $\{g(\cdot - k) \, | \, k \in \Z \}$ is an orthonormal system, then the measure of $supp(\hat{g}) \geq 1$. Equality holds iff $|\hat{g}| = \chi_K$ on some measurable set $K \subseteq \R$ of measure 1.
\ec}

\bp
By Plancherel's Theorem, a result in Fourier analysis, $||g||_2 = 1$ implies $|| \hat g ||_2 = 1$. From Lemma 1, we know that $|\hat g(x)| \leq 1$ almost everywhere on $\R$. Thus 
\bq
m(supp \, \hat g) = \int_{supp \, \hat g} 1 \, d\xi \geq \int_{\R} \hat g(\xi) \, d\xi = 1.
\eq 
Suppose $m(supp \, \hat g)=1$, and $|\hat g(\xi)|<1$ on a set $E$ of positive measure. Then we have
\begin{align*}
1 &= \int_{supp \, \hat g} |\hat g(\xi)|^2 \, d\xi = \int_E |\hat g(\xi)|^2 \, d\xi + \int_{supp \, \hat g - E} |\hat g(\xi)|^2 \, d\xi \\
&< m(E)+m(supp \, \hat g - E)=m(supp \, \hat g) = 1,
\end{align*}
a contradiction. Thus $|\hat g(\xi)|=\chi_K$, where $K=supp \, \hat g$, $m(K)=1$.
\ep

\section{Decomposition of $L^2(\R)$}

We will now construct an orthonormal wavelet from an MRA. Let $W_0$ be the orthogonal complement of $V_0$ in $V_1$: $V_1 = V_0 \oplus W_0$. If we dilate $W_0$ by $2^j$, we obtain the subspaces $W_j$ of $V_{j+1}$ such that $V_{j+1}=V_j \oplus W_j$ for all $j \in \Z$. Since $V_j \to 0$ as $j \to -\infty$, we see that

\bq \tag{3} V_{j+1}=V_j \oplus W_j = \bigoplus_{l=-\infty}^j W_l \quad \forall \, j \in \Z. \eq

Since $V_j \to L^2(\R)$ as $j \to +\infty$, we also have

\bq \tag{4} L^2(\R)=\bigoplus_{j=-\infty}^{\infty}W_j. \eq

We want to find $\psi \in W_0$ such that $\{ \psi(\cdot - k): \, k \in \Z \}$ is an orthonormal basis of $W_0$, then $\{2^{j/2} \psi(2^j \cdot - k): \, k \in \Z \}$ is an orthonormal basis for $W_j$, for all $j \in \Z$. Then it is an orthonormal wavelet basis for $L^2(\R)$, by (4).

Consider $V_0 = W_{-1} \oplus V_{-1}$, and observe that $\frac{1}{\sqrt{2}} \varphi(\frac{\cdot}{2}) \in V_{-1} \subseteq V_0$. We can express this function in terms of the basis $\{\varphi(\cdot - k): \, k \in \Z \}$: $\frac{1}{\sqrt{2}} \varphi(\frac{x}{2})=\sum_{k \in \Z} \alpha_k \varphi(x - k)$ where $(\sum_{k \in \Z} |\alpha_k|^2)^{1/2} < \infty$ and the convergence is in $L^2[0,1]$.

Taking Fourier Transforms,

\bq \tag{5} \hat{\varphi}(2\xi)=\hat{\varphi}(\xi)\sum_{k \in \Z} \alpha_k e^{2\pi ik \xi} = \hat{\varphi}(\xi)m_0(\xi), \eq

where $m_0(\xi)$ is the \textbf{low-pass filter} associated with the scaling function $\varphi$.

\section{Characterization of $V_0$ and $V_{-1}$}

We continue the construction of $\psi$. It is natural to apply Lemma 1 to the scaling function $\varphi$; then $\sum_{k \in \Z} |\hat \varphi(2\xi + k)|^2 = 1$ a.e., implying, by (5), $\sum_{k \in \Z} |\hat \varphi(\xi + k/2)|^2 |m_0(\xi + k/2)|^2 = 1$ a.e.. Taking the sum of the LHS separately over the even and odd integers, we have

\bq \tag{6} 1 = |m_0(\xi)|^2 \sum_{l \in \Z} |\hat \varphi(\xi + l)|^2 + |m_0(\xi + 1/2)|^2 \sum_{l \in \Z}|\hat \varphi(\xi + l + 1/2)|^2 = |m_0(\xi)|^2 + |m_0(\xi + 1/2)|^2.
\eq

This is known as the \textbf{Smith-Barnwell equality}. ($\overline{m_0(\xi + 1/2)}$ is known as the \textbf{high-pass filter}.)

If $f \in V_{-1}$, then $f(x)=\frac{1}{\sqrt{2}}\sum_{k \in \Z}c_k \varphi(x/2-k)$. Hence,

\begin{align*} 
\hat f(\xi)  &= \frac{1}{\sqrt{2}} \sum_{k \in \Z} c_k \int_{\R} \varphi(x/2-k)e^{-2\pi i x \xi} \, dx \\
&= \sum_{k \in \Z} c_k \int_{\R} \varphi(x/2-k)e^{-2\pi i(x/2-k)2\xi} \, d(x/2-k) \cdot \sqrt{2}e^{-2 \pi ik(2\xi)} \\
&= \sqrt{2}\sum_{k \in \Z} c_k \int_{\R} \varphi(y)e^{-2 \pi iy(2\xi)} \, dy \cdot e^{-2\pi ik(2\xi)} \\
&= \sqrt{2} \hat \varphi(2\xi)\sum_{k \in \Z}c_k e^{-2 \pi ik (2\xi)} \quad \text{ (by (5)) } \\
&= m(2\xi)m_0(\xi) \hat \varphi(\xi),
\end{align*}
where $m(2\xi)=\sqrt{2} \sum_{k \in \Z} c_k e^{-2 \pi i k (2\xi)}$.

We thus have a characterization of $V_0$,
\bq \tag{7} V_0 = \{f \in L^2(\R): \hat f(\xi)=l(\xi)\hat \varphi(\xi) \text{ for some 1-periodic } l \in L^2[0,1]\}, \eq

and of $V_{-1}$:
\bq \tag{8} V_{-1} = \{f \in L^2(\R): \hat f(\xi)=m(2\xi)m_0(\xi)\hat \varphi(\xi) \text{ for some 1-periodic } m \in L^2[0,1]\}. \eq

\section{Characterization of $W_{-1}$ and $W_0$}

We continue with the construction of the wavelet $\psi$. The elements of $W_{-1}$ are those $f \in V_0$ that are orthogonal to $V_{-1}$. Let $u:V_0 \to L^2[0,1]$ be defined by $u(f)=l$, $||u(f)||_{L^2[0,1]}^2=||l||_{L^2[0,1]}^2 = \sum_{k \in \Z} |d_k|^2$. If $f$ is perpendicular to $V_{-1}$, we must have that $l(\xi)$ is orthogonal to $m(2\xi)m_0(\xi)$ for all one-periodic $m \in L^2(\R)$. Then:

\bq 
0 = \int_0^1 l(\xi)\overline{m(2\xi)m_0(\xi)} \, d\xi = \int_0^{1/2} \overline{m(2\xi)}[l(\xi)\overline{m_0(\xi)}+l(\xi+1/2)\overline{m_0(\xi+1/2)}] \, d\xi.
\eq

The above equation implies that the 1-periodic function in the square brackets is orthogonal to all 1-periodic square integrable functions; that is, $l(\xi)\overline{m_0(\xi)}+l(\xi+1/2)\overline{m_0(\xi+1/2)}=0$ for almost every $\xi \in \mathbb{T}$. Hence, we must have
\bq \tag{9}
(l(\xi),l(\xi+1/2)) = -\lambda(\xi + 1/2)(\ol{m_0(\xi+1/2)},-\ol{m_0(\xi)}),
\eq
for a.e. $\xi$ and an appropriate $\lambda(\xi)$. We perform a change of variables: let $\xi \mapsto \xi + 1/2$. Then
\bq
(l(\xi+1/2),l(\xi)) = -\lambda(\xi + 1)(\ol{m_0(\xi)},-\ol{m_0(\xi+1/2)}),
\eq
by the 1-periodicity of $m_0$ and $l$. But the equality is equivalent to
\bq \tag{10}
(l(\xi),l(\xi+1/2)) = \lambda(\xi + 1)(\ol{m_0(\xi+1/2)},-\ol{m_0(\xi}))
\eq
for a.e. $\xi$, by a simple change of basis and adjusting for a factor of $-1$ on the RHS.

From (6), the Smith-Barnwell Equality, we know that the vector
\bq
(\ol{m_0(\xi+1/2)},-\ol{m_0(\xi)})
\eq
has norm $1$ for a.e. $\xi$. Combined with equations (9) and (10), we have $\lambda(\xi)=-\lambda(\xi + 1/2)$. Hence $\lambda$ is 1-periodic on $L^2[0,1]$, so $\exists$ a 1-periodic $s \in L^2[0,1]$ such that $\lambda(\xi)=e^{2\pi i\xi}s(2\xi)$. Rewriting, we get $s(\xi) = e^{-2\pi i\xi/2}\lambda(\xi/2)$. Then we obtain

\bq \tag{11}
l(\xi)=e^{2\pi i\xi}s(2\xi)\overline{m_0(\xi+1/2)}.
\eq

This gives us a characterization of $W_{-1}$:

\bq
W_{-1}=\{f : \hat f(\xi) = e^{2\pi i \xi} s(2\xi) \ol{m_0(\xi+1/2)} \hat \varphi(\xi) \textrm{ for a 1-periodic } s \in L^2[0,1] \},
\eq

which, in turn, characterizes $W_0$:

\blockg{\bl[12]
If $\varphi$ is a scaling function for an MRA $(V_j)_{j \in \Z}$, and $m_0$ is the associated low-pass filter, then
\bq
W_{0}=\{f : \hat f(2\xi) = e^{2\pi i \xi} s(2\xi) \ol{m_0(\xi+1/2)} \hat \varphi(\xi) \textrm{ for a 1-periodic } s \in L^2[0,1] \}.
\eq
\el}

\section{Characterization of orthonormal wavelets in $W_0$}

We are almost done with the construction of $\psi$. In Lemma 12, if we take $s \equiv 1$, that is, 
\bq \tag{13}
\hat \psi(2\xi) = e^{2\pi i \xi} \ol{m_0(\xi + 1/2)} \hat \varphi(\xi),
\eq
we claim we have found the desired orthonormal wavelet. In fact, we have completely characterized the orthonormal wavelets in $W_0$:

\blockg{\bprop[14]
Suppose $\varphi$ is a scaling function for an MRA $(V_j)_{j \in \Z}$, and $m_0$ is the associated low-pass filter, then a function $\psi \in W_0 = V_1 \cap V_0^{\perp}$ is an orthonormal wavelet for $L^2(\R)$ if and only if
\bq
\hat \psi(2\xi) = e^{2\pi i \xi} \nu(2\xi) \ol{m_0(\xi + 1/2)} \hat \varphi(\xi)
\eq
a.e. on $\R$, for some 1-periodic, measurable, a.e. unimodular function $\nu \in L^2[0,1]$.
\eprop}

\bp
Clearly $\psi \in W_0$, since we assumed that $\nu \in L^2[0,1]$. For any $g \in W_0$, by our characterization of $W_0$, $\exists \, s \in L^2[0,1])$, one-periodic, such that $\hat g(2\xi) = e^{2 \pi i \xi} s(2\xi)\ol{m_0(\xi+1/2)}\hat \varphi(\xi)$. This gives us
\bq
\hat g(\xi)=\frac{s(\xi)}{\nu(\xi)}e^{2 \pi i \xi/2} \nu(\xi) \ol{m_0(\xi/2+1/2)}\hat \varphi(\xi/2) = \frac{s(\xi)}{\nu(\xi)} \hat \psi(\xi) = s(\xi)\ol{\nu(\xi)}\hat \psi(\xi).
\eq 
Since $s\ol{\nu} \in L^2[0,1]$, we can write $s(\xi)\ol{\nu(\xi)}=\sum_{k \in \Z} c_k e^{-2\pi ik \xi}$ for a sequence $(c_k)_{k \in \Z} \in l^2(\Z)$, and obtain
\bq
g(x)=\sum_{k \in \Z}c_k \psi(x-k),
\eq
proving that $\{\psi(\cdot - k): k \in \Z \}$ generates $W_0$. To prove that this system is orthonormal, we show that $\hat \psi$ satisfies the equality in Lemma 1:
\begin{align*}
\sum_{k \in \Z} |\hat \psi(\xi + k)|^2 &= \sum_{k \in \Z} |\hat \varphi(\xi/2+k/2)|^2|m_0(\xi/2+ k/2 + 1/2)|^2 \\
&= \sum_{k \in \Z} |\hat \varphi(\xi/2+l)|^2|m_0(\xi/2+ l + 1/2)|^2 \\
& \quad + \sum_{k \in \Z} |\hat \varphi(\xi/2+l+1/2)|^2|m_0(\xi/2+ l + 1)|^2 \\
&= |m_0(\xi/2+1/2)|^2+|m_0(\xi/2)|^2=1,
\end{align*}
by summing over even and odd integers separately, using the 1-periodicity of $m_0$, Lemma 1 for $\varphi$ and (6), the Smith-Barnwell Equality, for $m_0$.

We have already observed that if $\{\psi(\cdot - k): k \in \Z \}$ is an orthonormal basis for $W_0$, then $\{2^{j/2}\psi(2^j\cdot - k): k \in \Z \}$ is an orthonormal basis for $W_j$. Hence (4) shows that $\psi$ is an orthonormal wavelet for $L^2(\R)$, as desired.

Now to show that all $\psi \in W_0$ are described by (13). Take $\psi \in W_0$. Then by Proposition 14, $\exists$ a $1$-periodic function $\nu \in L^2[0,1]$, such that
\bq
\hat \psi(\xi)=e^{2 \pi i \xi/2} \nu(\xi) \ol{m_0(\xi/2+1/2)}\hat \varphi(\xi/2).
\eq
If $\psi$ is an orthonormal wavelet, then the orthonormality of $\{\psi(\cdot - k): k \in \Z \}$ gives us
\begin{align*}
1 &= \sum_{k \in \Z} |\hat \psi(\xi+k)|^2 = \sum_{k \in \Z} |\nu(\xi)|^2 |m_0(\xi/2+k/2+1/2)|^2 |\hat \varphi(\xi/2+k/2)|^2 \\
&= |\nu(\xi)|^2 (\sum_{l \in \Z} |m_0(\xi/2+1/2)|^2 |\hat \varphi(\xi/2+l)|^2 + \sum_{l \in \Z} |m_0(\xi/2)|^2 |\hat \varphi(\xi/2+l+1/2)|^2) \\
&= |\nu(\xi)|^2(|m_0(\xi/2+1/2)|^2 + |m_0(\xi/2)|^2) = |\nu(\xi)|^2 \quad \text{for a.e. } \xi \in [0,1],
\end{align*}
and we are done.
\ep

\end{document}